\documentclass[12pt, letterpaper]{article}
\usepackage{fullpage}
\usepackage{amsmath, amssymb,amsthm}
\usepackage{enumerate}

\newtheorem{theorem}{Theorem}[section]
\newtheorem{lemma}[theorem]{Lemma}

\numberwithin{equation}{section}
 
\newenvironment{definition}[1][Definition]{\begin{trivlist}
\item[\hskip \labelsep {\bfseries #1}]}{\end{trivlist}}

\newenvironment{remark}[1][Remark]{\begin{trivlist}
\item[\hskip \labelsep {\bfseries #1}]}{\end{trivlist}}

\def\Z{{\mathbb Z}}         % integers
\def\R{{\mathbb R}}       % reals
\def\rect{{\mathcal R}}

\title{Maximal Functions for Lacunary Dilation Structures}
\author{Patrick LaVictoire}
\begin{document}

\maketitle
\begin{abstract}
If $\mu$ is a smooth density on a hypersurface in $\R^n$ whose curvature never vanishes to infinite order, and $A\in GL_d(\R)$ is a matrix whose eigenvalues all have absolute value greater than 1, let $\mu_k$ be the dilate of $\mu$ by $A^k$. We prove that $Tf=\sup_k f\ast \mu_k$ is bounded from a corresponding version of $H^1$ to weak $L^1$. %We also prove examples of a similar nature for measures on surfaces of higher codimension.
\end{abstract}

Consider $A \in GL_d(\R)$ whose eigenvalues all have absolute value greater than 1; this defines a dilation structure on $\R^d$. Let $\mu$ be a smooth density supported on a compact hypersurface $M\subset\R^d$ whose Gaussian curvature never vanishes to infinite order (by scaling, we may assume that this hypersurface is contained in the unit ball), and define the measures $\mu_k$ by
\begin{eqnarray}
\langle \mu_k, f\rangle = \langle \mu, f\circ A^k \rangle.
\end{eqnarray}
We will consider the maximal function
\begin{eqnarray}
{\cal M}f(x)=\sup_{k \in \Z}|\mu_k\ast f(x)|.
\end{eqnarray}
Littlewood-Paley theory quickly shows that $\cal M$ is bounded on $L^p$ for $1<p<\infty$, but the endpoint behavior has been a topic of some interest. Extrapolation arguments show that ${\cal M}:L\log L\to L^{1,\infty}$, but go no further. In a seminal paper, Christ \cite{MC} gave an endpoint result in two special cases: if $\mu$ is the surface measure on the sphere $S^{d-1}\subset \R^d$ and $A=2I$, then $\cal M$ is bounded from the real Hardy space $H^1(\R^d)$ to $L^{1,\infty}$; and if $M$ is the parabola $\{(t,t^2):1\leq t\leq 2\}\subset\R^2$ and $A(x,y)=(2x,4y)$, then $\cal M$ is bounded from the parabolic Hardy space $H^1_p(\R^d)$ (defined in terms of the parabolic dilation structure) to $L^{1,\infty}$. %The second result, in particular, involves a few ingredients that will be necessary for our result: a stopping-time construction after the fashion of the Calder\'{o}n-Zygmund decomposition in the parabolic Hardy space, but in such a way that it respects features of the Euclidean metric as well. Then, the curvature of the
\\
\\ It is relatively simple to extend the first of these results; if $A=2I$ and the Gaussian curvature of $M$ never vanishes, then ${\cal M}:H^1\to L^{1,\infty}$. As for the second, Seeger, Tao and Wright \cite{STW} obtained, for nearly our general hypotheses, a different endpoint space: namely that ${\cal M}:L\log \log L\to L^{1,\infty}$. The construction localizes the measure in order to handle vanishing curvature of $M$.
\\
\\Heo \cite{Heo} found room in the original construction of Christ to localize the measure $\mu$, showing for the usual isotropic dilation structure $A=2I$ that $\cal M$ is bounded from the Hardy space $H^1(\R^d)$ to weak $L^1$ (so long as the Gaussian curvature of $\mu$ never vanishes to infinite order); this was further generalized by Seeger and Wright \cite{SW}. The author has further extended this result with the full strength of the stopping-time argument, in order to obtain the following result:
\begin{theorem}
Let $M$, $\mu$, and $A$ be as above, and let $H_A^1(\R^d)$ be the Hardy space defined in terms of the dilation matrix $A$. Then $\cal M$ is bounded from $H_A^1(\R^d)$ to $L^{1,\infty}(\R^d)$.
\end{theorem}
The anisotropic Hardy space $H_A^1$ is best defined using the dyadic decomposition for spaces of homogeneous type, as in \cite{ChristTB}, Theorem 11. (Since all eigenvalues have norm greater than 1, $\rho(x,y)=\exp(\inf\{k: y-x\in A^k\})$ is a quasi-metric.)
\\
\\By the properties of that decomposition, there are constants $c$ and $C$ such that for any dyadic `cube' there exists an integer $k$ such that the `cube' contains a translate of $cA^kB_1(0)$ and is contained in a translate of $CA^kB_1(0)$. Thus we shall equivalently define this Hardy space in terms of a fixed grid which is invariant under $A$.
\begin{definition} 
For each $k\in\Z$ and $\vec n\in \Z^d$, let $Q^k_{\vec n}=A^k\left([0,1]^d+\vec n\right)$, and call $a_Q$ an atom if $a_Q$ is supported on such a set $Q$ and $\|a_Q\|_\infty\leq|Q|^{-1}=|\det A|^{-k}$; then $H_A^1$ is the completion of the finite sums $f=\sum_Q \lambda_Qa_Q$ in the norm $\|f\|_{H_A^1}=\sum_Q|\lambda_Q|$.
\end{definition}
(Technically, as in \cite{MC} we use two such grids, the second translated in space, to decompose a finite sum of the ``original'' atoms into finite sums with comparable norms in the grids.)
\\
\\We will need to keep track of two facts about our dilation structure: the volume and the diameter of our cubes.
\begin{definition}
Let $a:=|\det A|$; then the volume of the set $A^\tau([0,1]^d)$ is equal to $a^\tau$ for each $\tau\in\Z$.
\end{definition}
\begin{definition}
Let $r$ denote the minimum of the absolute values of the eigenvalues of $A$, and let $n$ denote the size of the largest block in the Jordan decomposition of $A$ whose eigenvalue has absolute value $r$. Then for $\tau\leq0$, the diameter of the set $A^\tau([0,1]^d)$ is comparable to $r^\tau |\tau|^n$.
\end{definition}

\section{Stopping Time Construction}

Here we will prove a more general version of the stopping-time lemma used in \cite{MC} to construct the exceptional set for the decomposition.
\\
\\ Now if ${\mathcal D}_\sigma$ is the dyadic grid of sidelength $2^\sigma$, we would like to define the grid $\rect_{\sigma,\tau}$ to be the image of ${\mathcal D}_\sigma$ under $A^\tau$. (Actually, we will replace $A$ with a constant integer multiple $A^n$, such that $A^n$ maps $B(0,1)$ into $B(0,\frac12)$. This requires us to divide our maximal function $\cal M$ into $n$ pieces, which of course does not affect our result.) We further let $\rect$ denote the union of all $\rect_{\sigma,\tau}$, and $\rect_0$ denote the strictly anisotropically dilated cubes $\bigcup_\tau \rect_{0,\tau}$ (the supports of our atoms $a_Q$).
\\
\\ Note that we have defined $\sigma$ differently from \cite{MC} in this more general context. Rather than designating the length of the side parallel to the $x$-axis, we have let it denote the number of isotropic dilations that separate $q$ from an element of $\rect_0$. In what follows, we will only need to consider $\sigma\leq0$.
\\
\\ As in \cite{MC}, we will build an exceptional set out of tendrils as well as other components. The first step is the standard anisotropic Whitney decomposition:
\begin{lemma}
\label{streets}
For any $\alpha>0$ and any collection of cubes $Q$ belonging to $\rect_0$, with associated positive scalars $\lambda_Q$, there exists a collection of pairwise disjoint cubes ${\cal S}\subset\rect_0$ such that
\begin{enumerate}
\item \label{turned} $\displaystyle \sum_{Q\subset S^*} \lambda_Q \leq C\alpha|S|$ for all $S\in{\cal S}$
\item $\sum |S|\leq \alpha^{-1}\sum \lambda_Q$
\item \label{others}$\displaystyle \left\| \sum_{Q\not\subset S^* \forall S}\lambda_Q\frac{\chi_Q}{|Q|}\right\|_\infty\leq \alpha$.
\end{enumerate}
\end{lemma}
The proof of this lemma is mostly standard, though we must use $S^*$ in places because we cannot count on the `children' of $S$ to be nested in $S$. We can safely discard the $Q$ not belonging to any of these $S^*$, by condition \ref{others}; although it is not true that each remaining $Q$ belongs to $S^*$ for a unique $S\in{\cal S}$, each $S$ whose double contains $Q$ will have identical dimensions. (Technically speaking, this is a consequence of the standard proof of this theorem.)
\\
\\ We define the tendril $T(q)$ of any $q\in\rect_{\sigma,\tau}$ by
\begin{eqnarray}
T(q):= q^{**}+\bigcup_{k\leq \tau(q)+2}\text{supp }\mu_k
\end{eqnarray}
and note that $\chi_{q^*}\ast \mu_j$ is supported in $T(q)$. (Here $q^*$ is the expansion of $q$ by a factor of 2, and $q^{**}$ is the expansion of $q$ by a factor of 4.) Note as well that $T(q)$ is the image of $q'+\bigcup_{k\leq 2}\text{supp }\mu_k$ under $A^\tau$ for some $q'\in \rect_{\sigma,0}$, and that therefore $|T(q)|\lesssim 2^\sigma |A|^\tau$.
\\
\\ Now we will generalize the stopping-time lemma from \cite{MC}:
\begin{lemma}
\label{collaborate}
We are given $\alpha>0$, a finite collection $\cal S$ of pairwise disjoint $S\in\rect_0$, and a finite collection $\cal Q$ of $Q\in \rect_0$, such that each $Q\in\cal Q$ is contained in $S^*$ for some $S\in\cal S$. Corresponding to each $Q\in\cal Q$ we are also given $\lambda_Q>0$. Then there exist a measurable $E\subset\R^3$ and a function $\kappa:{\cal Q}\to\Z$ such that
\begin{enumerate}[i)]
\item \label{ice} $|E|\leq C\alpha^{-1}\sum \lambda_Q +C\sum |S|$
\item \label{baby} $\chi_Q \ast \mu_j$ is supported in $E$ for all $Q\in \cal Q$ and all $j<\kappa(Q)$
\item \label{vanilla} If $Q\subset S^*$, then $\kappa(Q)>\tau(S)$
\item \label{listen} For any $\tau\in \Z$ and $\sigma\leq0$ and $q \in \rect_{\sigma,\tau}$, \begin{eqnarray}\label{stop}
\sum_{Q\subset q^*: \kappa(Q)\leq \tau} \lambda_Q \leq C\alpha|T(q)|.
\end{eqnarray}
\end{enumerate}
\end{lemma}
\begin{proof}
The proof proceeds via a double induction on $\sigma$ and $\tau$. We will partition $\cal Q$ into two subcollections ${\cal C}_1$ and ${\cal C}_2$. To each $Q\in{\cal C}_1$ we will associate a $q\in\cal R$ with $Q\subset q^*$, and we will (initially) set $\kappa(Q)=\tau(q)+1$. Each $Q\in {\cal C}_2$ will have $\kappa(Q)=\tau(S)+1$, where $S^*$ contains $Q$. (As mentioned above, this uniquely determines $\tau(S)$.)
\\
\\Select $\tau_0$ larger than $\tau(Q)$ for all $Q\in\cal Q$, such that $\alpha |A|^{\tau}>\displaystyle\sum_{Q\in\cal Q} \lambda_Q$. Initialize $\tau=\tau_0-1$ and $\sigma=0$, and define ${\cal Q}(0,\tau_0)={\cal Q}$; we will define collections ${\cal Q}(\sigma, \tau)$ by removing elements at each step so that ${\cal Q}(\sigma,\tau)\subset{\cal Q}(\sigma+1, \tau)\subset{\cal Q}(\sigma', \tau+1)$ for any $\tau\leq\tau_0$ and any $\sigma$ and $\sigma'\leq0$. Furthermore, we will define $\Lambda_{\sigma, \tau}(q)=\displaystyle\sum_{\substack{Q\subset q^* \\ Q \in {\cal Q}(\sigma,\tau)}}\lambda_Q$ for each $q\in\rect_{\sigma,\tau}$.
\\
\\For each $\tau$ in descending order from $\tau_0$, proceed by descent on $\sigma$; for each fixed $\sigma$, select all $q\in \rect_{\sigma,\tau}$ such that $\Lambda_{\sigma,\tau}(q)>\alpha 2^{\sigma}|A|^{\tau}$. Any $Q$ contained in $q^*$ for a selected $q$ is classified into ${\cal C}_1$ and assigned to one of the selected $q$; $\kappa(Q)$ is defined to be $\tau(q)+1$. (Again, the assignment of $q$ will not be unique, but the assignment of $\kappa(Q)$ is.) Then ${\cal Q}(\sigma-1,\tau)$ is defined to consist of all $Q\in{\cal Q}(\sigma,\tau)$ which were not classified at this step.
\\
\\ Eventually, $-\sigma$ is so large that no $Q\in \cal Q$ can be contained in $q^*$ for any $q\in \rect_{\sigma,\tau}$. When this happens, all unassigned $Q$ with the dimensions $\tau(Q)=\tau$ are classified into ${\cal C}_2$, assigned to some $S$ with $Q\subset S^*$, and given $\kappa(Q)=\tau(S)+1$. Then $\tau$ is incremented down by 1, ${\cal Q}(0,\tau-1)$ is defined to consist of the remaining unclassified $Q$ (i.e. all $Q\in{\cal Q}(-\infty,\tau)$ such that $\tau(Q)<\tau$), and we start descending in $\sigma$ again. This process repeats until we reach $\tau$ smaller than $\tau(Q)$ for all $Q$, at which point all of $\cal Q$ has been classified.
\\
\\ Throughout this process, we have ensured the usual stopping-time condition $\Lambda_{\sigma,\tau}(q)\leq C\alpha 2^{\sigma}|A|^{\tau}\leq C\alpha|T(q)|$ for all $q\in\rect$, since otherwise a parent of $q$ would have been chosen instead. (This is true as well if $\sigma=0$, since in that case we may consider the anisotropic parents of $q$.) Of course, the left-hand side of (\ref{stop}) is precisely $\Lambda_{\sigma,\tau}(q)$, so condition (\ref{listen}) is verified.
\\
\\Now we define $E_2$ to be the union of all the tendrils $T(q)$ for all $q$ selected in the process and $E_1$ to be the union of all the quadruples $S^{**}$ of $S\in\cal S$. We already know that $|E_2|$ is appropriately bounded in size. For $|E_1|$, each $Q$ is assigned to one of at most $3^d$ cubes $q$. Therefore, summing over the selected $q$,
\begin{eqnarray*}
|E_1|\leq \sum_q |T(q)| \lesssim \sum_q \alpha^{-1}\Lambda_{\sigma(q),\tau(q)}(q)\lesssim  \alpha^{-1}\sum_{Q\in{\cal C}_1}  \lambda_Q.
\end{eqnarray*}
Therefore, with $E=E_1\cup E_2$, we have satisfied condition (\ref{ice}). Condition (\ref{baby}) is clearly true for $Q\in {\cal C}_1$, and since $\mu$ is supported in the unit ball, it is trivial to show for $Q\in{\cal C}_2$ as well. Unfortunately, condition (\ref{vanilla}) need not hold; however, we can repair this by replacing the current $\kappa$ with $\max(\kappa(Q),\tau(S)+1)$. This does not affect \ref{ice}); (\ref{baby}) is still true because $\kappa(Q)$ is either unchanged or replaced with $\tau(S)+1$, which we noted is fine; it makes (\ref{vanilla}) trivially true; and it preserves (\ref{listen}) because fewer $Q$ will now be summed over on the left. Thus we are done.
\end{proof}

\section{Proof of Main Theorem}

We will use the following lemma from Seeger, Tao and Wright (\cite{STW}, Lemma 2.5):

\begin{lemma}
\label{bowie}
Let $\psi\in C^\infty([-1,1]^{d-1})$ be real-valued, with $\sup_{|\alpha|\leq 3}|\partial^\alpha \psi(x)|\leq A\leq 1$ on $[-1,1]^{d-1}$. Suppose $|\det \psi''(y_0)|\geq \beta$, and $Q\subset[-1,1]^{d-1}$ is a cube of sidelength $\epsilon_1 \beta$ containing $y_0$, where $\epsilon_1\leq[10(d-1)^4A]^{-1}$.
Let $\chi\in C^\infty(Q)$ with $\|\partial^\alpha \chi\|_\infty \leq c_\alpha(\epsilon_1 \beta)^{-|\alpha|}$. Define the measure $\nu$ on $\R^d$ by
\begin{eqnarray*}
\langle \nu, f\rangle = \int \chi(y) f(y,\psi(y)) dy
\end{eqnarray*}
and define its reflection by $\langle \tilde\nu, f\rangle=\langle \nu, f(-\cdot)\rangle$.
Then there are constants $C_\alpha$ so that
\begin{eqnarray}
|\partial_x^\alpha[\nu\ast\tilde\nu](x)|\leq C_\alpha \beta^{d-3-2|\alpha|}|x|^{-1-|\alpha|}.
\end{eqnarray}
\end{lemma}
As in \cite{MC}, we will use this regularity of the kernel convolved with its reflection to obtain an especially strong $L^2$ bound on a part of our operator. In \cite{Heo}, this idea was combined with that of partitioning $\mu$ into pieces of small support, and setting aside those which are ``bad'' in a certain sense; so long as we can bound the contribution of these ``bad'' pieces in $L^1$, we may assume quantitative conditions on the remaining pieces.
\\
\\ For each $s$, we use a smooth partition of unity to write $\mu=\sum_{\rho\in I_s}\mu_\rho^s$, where $\mu_\rho^s$ is supported on a ball $B^s_\rho$ of diameter $2^{-\epsilon s}$, and $|I_s|\lesssim 2^{(d-1)\epsilon s}$. (We take $0<\epsilon\ll \log(r)$, which ensures that $B^s_\rho$ will have large diameter compared to the atoms it will be convolved with.)
\\
\\There are two types of pieces that we'd like to exclude. As in \cite{Heo}, we will set aside those pieces on which the curvature falls below $2^{-\epsilon s}$ (since the size of the curvature is used in Lemma \ref{bowie}). We will also need a certain ``transversality'' condition; the intersection of $Q\in\rect_{0,\tau}$ with a typical piece of $M$ should have measure comparable to the volume of $Q$ divided by the diameter of $Q$, and we will use this fact.
\begin{definition}
Let 
\begin{eqnarray}
I_s^1:=\{\rho: \min_{x\in\text{supp}(\mu_\rho^s)}|K(x)|<2^{-\epsilon s}\},
\end{eqnarray}
where $K$ is the Gaussian curvature of the manifold, and for $0<\zeta\ll\epsilon$ let 
\begin{eqnarray}
I_s^2:=\left\{\rho: \exists Q\in\rect_0 \text{ such that }\mu^s_\rho(Q) > 2^{\zeta s}|Q|\text{diam}(Q)^{-1}\right\}.
\end{eqnarray}
\end{definition}
Then we claim that the contribution of these sets in $L^1$ will be summable:
\begin{lemma}
\label{mercury}
There exists $\eta>0$ such that $|I_s^1\cup I_s^2|\lesssim 2^{((d-1)\epsilon-\eta)s}$.
\end{lemma}
\begin{proof}
The bound on $|I_s^1|$ follows from the nonvanishing curvature of $M$, as shown in \cite{Heo}.
\\
\\For the bound on $|I_s^2|$, we first identify the direction of slowest contraction under the dilation group, noting that Jordan blocks can induce a logarithmic factor. Consider the real Jordan form of $A$; among the blocks whose eigenvalues have norm $r$, choose the one whose block size $n$ is maximal. (A $2n\times 2n$ complex Jordan block counts as size $n$ here.)
\\ 
\\In the corresponding eigenspace, there exists a unit vector $\vec v$ and a subspace $W$ (of dimension 1 or 2, depending on whether it is a real or complex eigenvalue) such that $A^\tau\vec v\asymp r^\tau|\tau|^n$ and the distance from $r^{-\tau}|\tau|^{-n}A^\tau\vec v$ to $W$ tends to 0 as $\tau\to-\infty$.
\\
\\ Now if $\vec N_x$ is the normal vector to $M$ at $x$, and $|\langle\vec N_x,\vec v\rangle|\geq c>0$ for all $x\in M\cap B_\rho$, then $\mu^s_\rho(M\cap Q)\lesssim c^{1-d}|\det A|^\tau |A^\tau\vec v|^{-1}$. The assumption on the curvature ensures that the estimate follows.
\end{proof}
Now let $f\in H^1(\R^d)$ be a finite sum $f=\sum_Q \lambda_Q a_Q$, where the $Q$ are all elements of $\rect_0$, with $\lambda_Q>0$, $\|a_Q\|_\infty\leq|Q|^{-1}$, and $\int a_Q =0$ for all $Q$. We want to show
\begin{eqnarray*}
|\{x: {\cal M}f(x)>2\alpha\}|\lesssim \alpha^{-1}\sum_Q \lambda_Q.
\end{eqnarray*}
We apply Lemma \ref{streets} to obtain the collection $\cal S \subset \rect_0$, and note that if $g$ denotes the sum of all $\lambda_Qa_Q$ for $Q$ not contained in any $S^*$, then by (\ref{others}), $\|g\|_\infty\leq\alpha$ and thus $\|{\cal M}g\|_\infty\leq\alpha$. Thus we may assume that all of the $Q$ are contained in $S^*$ for some $S\in\cal S$.
\\
\\ We now apply Lemma \ref{collaborate} to our collection $\cal Q$ of such cubes and our collection $\cal S$, obtaining the exceptional set $E$ and the function $\kappa:{\cal Q}\to\Z$. Since $|E|\lesssim \alpha^{-1}\sum_Q \lambda_Q$, it suffices to prove that
\begin{eqnarray*}
|\{x\notin E: \sup_k |\mu_k\ast(\sum_Q \lambda_Q a_Q)(x)|>\alpha\}|\lesssim \alpha^{-1}\sum_Q \lambda_Q.
\end{eqnarray*}
By (\ref{baby}), we see that
\begin{eqnarray*}
|\{x\notin E: \sup_k |\mu_k\ast(\sum_Q \lambda_Q a_Q)(x)|>\alpha\}|&\leq& |\{x: \sup_j |\mu_j\ast(\sum_{Q:\kappa(Q)\leq j} \lambda_Q a_Q)(x)|>\alpha\}|\\
&\leq& |\{x: \sup_j |\sum_{s=0}^\infty\mu_j\ast(\sum_{Q:\kappa(Q)= j-s} \lambda_Q a_Q)(x)|>\alpha\}|.
\end{eqnarray*}
We now partition $\mu$ as discussed above, depending on $s$; write
\begin{eqnarray*}
{\cal M}' f(x)&:=& \sup_j |\sum_{s=0}^\infty\sum_{\rho\in I_s^1\cup I_s^2}\mu_{\rho,j}^s\ast(\sum_{Q:\kappa(Q)= j-s} \lambda_Q a_Q)(x)|,\\
{\cal M}'' f(x)&:=& \sup_j |\sum_{s=0}^\infty\sum_{\rho\notin I_s^1\cup I_s^2}\mu_{\rho,j}^s\ast(\sum_{Q:\kappa(Q)= j-s} \lambda_Q a_Q)(x)|.
\end{eqnarray*}
Now by Lemma \ref{mercury}, we see that
\begin{eqnarray*}
|\{x: {\cal M}'f(x)>\alpha/2\}|&\leq& \frac2\alpha \|{\cal M}'f\|_{L^1}\leq \frac2\alpha \sum_j \sum_{s=0}^\infty\left\|\sum_{\rho\in I_s^1\cup I_s^2}\mu_{\rho,j}^s\right\|_1\left\|\sum_{Q:\kappa(Q)= j-s} \lambda_Q a_Q\right\|_1\\
&\lesssim& \alpha^{-1}\sum_{s=0}^\infty 2^{-\eta s}\sum_j\sum_{\kappa(Q)=j-s} \lambda_Q\lesssim\alpha^{-1}\sum_Q\lambda_Q.
\end{eqnarray*}
Thus we only need concern ourselves with ${\cal M}''$. By Chebyshev's inequality, it will suffice to prove that $\|{\cal M}''f\|_2^2\lesssim \alpha\sum_Q \lambda_Q$. Note that
\begin{eqnarray*}
|{\cal M}''f(x)|^2&\leq& \left(\sup_j \left|\sum_{s=0}^\infty\sum_{\rho\notin I_s^1\cup I_s^2}\left(\sum_{Q:\kappa(Q)= j-s} \lambda_Q a_Q\right)\ast\mu_{\rho,j}^s(x)\right|\right)^2\\
&\leq& \sum_j\left|\sum_{s=0}^\infty\sum_{\rho\notin I_s^1\cup I_s^2}\left(\sum_{Q:\kappa(Q)= j-s} \lambda_Q a_Q\right)\ast\mu_{\rho,j}^s(x)\right|^2
\end{eqnarray*}
and by Minkowski's inequality,
\begin{eqnarray*}
|{\cal M}''f(x)| &\leq& \sum_{s=0}^\infty\sum_{\rho\notin I_s^1\cup I_s^2}\left( \sum_j \left|\left(\sum_{Q:\kappa(Q)= j-s} \lambda_Q a_Q\right)\ast\mu_{\rho,j}^s(x)\right|^2\right)^{1/2}
\end{eqnarray*}
so that
\begin{eqnarray*}
\|{\cal M}''f\|_2 &\leq& \sum_{s=0}^\infty\sum_{\rho\notin I_s^1\cup I_s^2}\left( \sum_j \left\|\left(\sum_{Q:\kappa(Q)= j-s} \lambda_Q a_Q\right)\ast\mu_{\rho,j}^s\right\|_2^2\right)^{1/2}\\
&\lesssim& \sum_{s=0}^\infty 2^{(d-1)\epsilon s}\sup_{\rho\notin I_s^1\cup I_s^2}\left( \sum_j \left\|\left(\sum_{Q:\kappa(Q)= j-s} \lambda_Q a_Q\right)\ast\mu_{\rho,j}^s\right\|_2^2\right)^{1/2}.
\end{eqnarray*}
Therefore it suffices to prove that there exists $C<\infty$ and $\delta>0$ such that for any $s\geq 0$, $\rho\in I_s\setminus (I_s^1\cup I_s^2)$, and $j\in\Z$,
\begin{eqnarray}
\label{terror}
\left\|\left(\sum_{\kappa(Q)=j-s} \lambda_Q a_Q\right)\ast \mu_{\rho,j}^s\right\|_{L^2}^2 \leq C2^{-(2\epsilon(d-1)+\delta)s}\alpha \sum_{\kappa(Q)=j-s} \lambda_Q.
\end{eqnarray}
Finally, by scaling, we may assume $j=0$.

\subsection{Proof for $s=0$}

In order to handle our anisotropic dilation structure, as in \cite{MC} we must consider elements of $\rect_0$ as subsets of isotropic dyadic cubes, since the bounds we will obtain (by convolving the measure $\mu$ with its reflection) depend on the Euclidean distance.
\\
\\ For $s=0$, we will consider each $Q\in \rect_{0,\tau}$ with $\kappa(Q)=0$ as a subset of an isotropic cube $q\in\rect_{\sigma,0}$ such that $Q\subset q^*$ and diam$(Q)\approx r^\tau|\tau|^n\approx 2^\sigma$. Define $$A_q(x)=\displaystyle\sum_{\substack{Q\subset q^*, \, \kappa(Q)=0,\\ \text{diam}(Q)\approx 2^\sigma}} \lambda_Q a_Q(x) \text{ and } \lambda_q=\sum_{\substack{Q\subset q^*, \, \kappa(Q)=0\\ \text{diam}(Q)\approx 2^\sigma}} \lambda_Q.$$ (To avoid double-counting, we actually modify this a bit: we partition the cubes $Q\in \rect_{0,\tau}$ among the $q$ of the appropriate size.)
\begin{lemma}
\label{muppet}
For any $q\in \rect_{\sigma,0}$, $\|A_q \ast \mu_{\rho,0}^s\|_\infty\leq C 2^{-\sigma+\zeta s}\lambda_q$, and $\|A_q \ast \mu_{\rho,0}^s\|_1 \leq C2^{(\zeta +\epsilon(1-d))s}\lambda_q$.
\end{lemma}
\begin{proof}
Since $\|a_Q\|_\infty\leq |Q|^{-1}$, the definition of $I_s^2$ implies that for any $Q$ contained in $q$, \begin{eqnarray*}\|a_Q \ast \mu_{\rho,0}^s\|_\infty\lesssim |Q|^{-1}\mu_{\rho,0}^s(Q)\leq 2^{\zeta s}\text{diam }(Q)^{-1}=2^{-\sigma+\zeta s},\end{eqnarray*} which implies the first inequality. The second follows from the fact that $A_q \ast \mu_{\rho,0}^s$ is supported on a set of measure $\leq 2^{\sigma + \epsilon(1-d)s}$.
\end{proof}
This immediately implies a first estimate for any $q$ and $q'$ with $\sigma(q')\geq\sigma(q)$,
\begin{eqnarray}
\label{under}
\left|\langle A_q \ast \mu_{\rho,0}^0, A_{q'} \ast \mu_{\rho,0}^0\rangle\right|\leq C2^{-\sigma'}\lambda_q\lambda_{q'}.
\end{eqnarray}
We will obtain a second estimate on this quantity (which is stronger when $q$ and $q'$ are distant but weaker when they are near) using our assumption on the curvature of $M\cap B_\rho$ (from the definition of $I_s^1$) and Lemma \ref{bowie}.
\begin{lemma}
\label{sesame}
If $q$ and $q'$ are separated by a Euclidean distance of $d(q,q')\geq 2^{\sigma'}$, then
\begin{eqnarray}
\label{pressure}
\left|\langle A_q \ast \mu_{\rho,0}^s, A_{q'} \ast \mu_{\rho,0}^s\rangle\right|\leq C2^{\sigma'+\epsilon s(5-d)}d(q,q')^{-2}\lambda_q\lambda_{q'}.
\end{eqnarray}\end{lemma}
\begin{proof} By translation, we may assume that $q'$ is centered at 0.
\begin{eqnarray*}
\langle A_q \ast \mu_{\rho,0}^s, A_{q'} \ast \mu_{\rho,0}^s\rangle &=& \langle A_q , A_{q'} \ast \mu_{\rho,0}^s\ast \tilde\mu_{\rho,0}^s \rangle.
\end{eqnarray*}
Let $\varphi(x)=\mu_{\rho,0}^s\ast \tilde\mu_{\rho,0}^s(x)$; by Lemma \ref{bowie}, we see that $|\varphi(x)|\lesssim 2^{\epsilon s(3-d)}|x|^{-1}$ and $|\nabla \varphi(x)|\lesssim 2^{\epsilon s(5-d)}|x|^{-2}$. Now for any $x\in q$,
\begin{eqnarray*}
A_{q'}\ast \varphi(x) &=& \sum_{Q \text{ assigned to }q'}\lambda_Q\int a_Q(y)\varphi(x-y)\, dy\\
&=& \sum_{Q \text{ assigned to }q'}\lambda_Q\int a_Q(y)\left[\varphi(x-y)-\varphi(x)\right], dy\\
\end{eqnarray*}
by the cancellation of $a_Q$, and therefore
\begin{eqnarray*}
|A_{q'}\ast \varphi(x)|\lesssim \sum_{Q \text{ assigned to }q'}\lambda_Q \left(\text{diam}(Q)^{-1}\sup_{x\in q} |\nabla\varphi(x)|\right)\|a_Q\|_1 \lesssim \lambda_{q'} 2^{\sigma'+\epsilon s(5-d)} d(q,q')^{-2}
\end{eqnarray*}
which implies the result.
\end{proof}
Now
\begin{eqnarray*}
\left\|\left(\sum_{\kappa(Q)=0} \lambda_Q a_Q\right)\ast \mu_{\rho,0}^0\right\|_2^2 &=& \left\| \sum_q A_q \ast\mu_{\rho,0}^0\right\|_2^2\\
&\leq& 2\sum_{\sigma\leq\sigma'\leq0} \sum_{q\in \rect_{\sigma,0},q'\in\rect_{\sigma',0}}\left|\langle A_q \ast \mu_{\rho,0}^0, A_{q'} \ast \mu_{\rho,0}^0\rangle\right|\\
&=& \left(2\sum_{q'}\sum_{q\subset {q'}^{**}}\right)+\left(2\sum_{q'}\sum_{q\cap {q'}^{**}=\emptyset}\right) = I+II.
\end{eqnarray*}
For the sum $I$, we use the trivial bound (\ref{under}) and property (\ref{listen}) of Lemma \ref{collaborate} (note that $|T(q')|\approx 2^{\sigma'}$) to see
\begin{eqnarray*}
I&\leq& C\sum_{q'}\sum_{q\subset {q'}^{**}} 2^{-\sigma'}\lambda_{q'}\lambda_q\\
&\leq& C\sum_{q'}\alpha\lambda_{q'}\leq C\alpha\sum_{\kappa(Q)=0}\lambda_Q
\end{eqnarray*}
as desired. For the sum $II$, note that $2^{\sigma'}\leq d(q,q')\leq 2^3$ for all $q,q'$. Thus by (\ref{pressure}),
\begin{eqnarray*}
II &\leq& C\sum_{q'}\sum_{m=\sigma'}^3 \sum_{q: d(q,q')\approx 2^m} 2^{\sigma'}2^{-2m}\lambda_q\lambda_{q'}.
\end{eqnarray*}
Since the $q$ with $d(q,q')\approx 2^m$ are contained in one of $4^d$ isotropic cubes $q''\in\rect_{m,0}$ (or one of at most $32^d$ cubes in $\rect_{0,0}$, if $1<m\leq 3$), we apply property (\ref{listen}) to those isotropic cubes and see
\begin{eqnarray*}
II &\leq& C\sum_{q'}\sum_{m=\sigma'}^3 2^{\sigma'}2^{-2m}\alpha 2^{m}\lambda_{q'}\\
&\leq& C\alpha\sum_{q'}\lambda_{q'}\sum_{m=\sigma'}^3 2^{\sigma'-m}=C\alpha\sum_{\kappa(Q)=0}\lambda_Q
\end{eqnarray*}
and the proof is complete for $s=0$.

\subsection{Proof for $s>0$}

We must gain the factor of $2^{-(2\epsilon(d-1)+\delta)s}$ that appears in (\ref{terror}); we will do this by taking advantage of the better bounds in (\ref{stop}) of Lemma \ref{collaborate} for elements of $\rect_{\theta,-s}$, as compared to $\rect_{\sigma,0}$. We will therefore simply rewrite $q\in \rect_{\sigma,0}$ as a union of cubes in some $\rect_{\theta,-s}$ whenever it comes time to apply that bound. However, we run out of room to do this for pairs of cubes $(q,q')$ that are too far apart; in that case, we instead will use the Whitney decomposition (Lemma \ref{streets}).\\
\\Again, we assign $Q$ with $\kappa(Q)=-s$ to an isotropic  $q\in\rect_{\sigma,0}$ of approximately equal diameter, but we also assign it to an element of $\rect_{\theta,-s}$, where $\theta$ is the largest integer such that an element of $\rect_{\sigma,0}$ can contain an element of $\rect_{\theta,-s}$. That is to say, $2^{\theta}r^{-s}s^n\approx 2^\sigma$. (Note that $Q$ will indeed be contained in one of these, as it has been anisotropically dilated at least $s$ times.)
\\
\\ Since $\kappa(Q)=-s$ implies $\tau(Q)<-s$, the diameter of $Q$ is at most $r^{-s}s^{n}$, and thus $\theta\leq0$.
\\
\\ Now by (\ref{listen}), for $q' \in\rect_{\theta,-s}$,
\begin{eqnarray*}
\sum_{\substack{Q\subset q'^*, \, \kappa(Q)=-s,\\ \text{diam}(Q)\approx \text{diam}(q')}}\lambda_Q\leq C\alpha |T(q')|=C\alpha 2^{\theta}|A|^{-s}
\end{eqnarray*}
and there are $\lesssim 2^{d(\sigma-\theta)}|A|^s$ such $q'$ contained in $q\in \rect_{\sigma,0}$, so that
\begin{eqnarray}
\sum_{\substack{Q\subset q^*, \kappa(Q)=-s,\\ \text{diam}(Q)\approx \text{diam}(q)}}\lambda_Q\leq C\alpha 2^{d\sigma+(1-d)\theta}.
\end{eqnarray}
This gains a factor of $(r^{-s}s^n)^{d-1}$ over the direct estimate for $q\in\rect_{\sigma,0}$, which we can use to offset our losses due to decreased curvature or due to alignment with the direction of slowest contraction. First, we see by Lemma \ref{muppet} that
\begin{eqnarray*}
I &=& 2\sum_{q'}\sum_{q\subset {q'}^{**}}\left|\langle A_q \ast \mu_{\rho,0}^s, A_{q'} \ast \mu_{\rho,0}^s\rangle\right|\\
&\lesssim& \sum_{q'}\sum_{Q\subset {q'}^{**}, \kappa(Q)=-s}  2^{\zeta s -\epsilon(d-1)s}\lambda_Q 2^{-\sigma'+\zeta s}\lambda_{q'}\\
&\lesssim& \sum_{q'} \lambda_{q'} 2^{-\sigma' +2\zeta s -\epsilon(d-1)s}\sum_{Q\subset {q'}^{**}, \kappa(Q)=-s}\lambda_Q\\
&\lesssim& \sum_{q'} \lambda_{q'} 2^{-\sigma' +2\zeta s -\epsilon(d-1)s}\alpha 2^{\sigma'}(r^{-s}s^{-n})^{1-d}\\
&\lesssim& 2^{-(2\epsilon(d-1)+\delta)s}\alpha\sum_Q\lambda_Q,
\end{eqnarray*}
so long as $\epsilon\ll \log(r)$, $\zeta\ll\epsilon$ and $\delta\ll\epsilon$.
\\
\\This same gain will help us for the part of $II$ where we can write $q\in\rect_{m,0}$ as the union of elements of $\rect_{\theta,-s}$ with $\theta\leq0$. However, this leaves us with a number of terms that we cannot handle in this manner: those with $d(q,q')\geq r^{-s}s^n$. (We separate out these terms before putting the absolute values on inner products; just set aside all $q,q'$ with that much distance between them.) These terms we will deal with separately, aggregating them instead via the cubes $S$. %Why are these a problem, again?
\\
\\First, we will consider the pairs with $d(q,q')\leq r^{-s}s^n$; we call this $II'$, and by Lemma \ref{sesame} and the fact that $T(q'')\lesssim 2^m(r^s s^{-n})^{1-d}$ for $q''\in \rect_{\theta, -s}$ here, we see that
\begin{eqnarray*}
II'&\leq& \sum_{q'}\sum_{m=\sigma'}^{\log_2(r^{-s}s^n)} \sum_{q: d(q)\approx 2^m} \left|\langle A_q \ast \mu^s_{\rho,0}\ast \tilde\mu^s_{\rho,0}, A_{q'}\rangle\right|\\
&\lesssim& \sum_{q'}\sum_{m=\sigma'}^{\log_2(r^{-s}s^n)} \sum_{q: d(q)\approx 2^m} 2^{\sigma' +\epsilon s(5-d)} 2^{-2m}\lambda_q\lambda_{q'}\\
&\lesssim& \sum_{q'}2^{\epsilon s(5-d)} \lambda_{q'}\sum_{m=\sigma'}^{\log_2(r^{-s}s^n)}2^{\sigma'-2m}\cdot\alpha 2^m(r^{-s}s^n)^{d-1}\\
&\lesssim& \alpha (r^{-s}s^n)^d 2^{\epsilon(5-d)s}\sum_{q'}\lambda_{q'}\\
&\lesssim& 2^{-(2\epsilon(d-1)+\delta)s}\alpha\sum_Q\lambda_Q.
\end{eqnarray*}
(Note that we used $\epsilon\ll \log(r)$ in the final step.)
\\
\\We still have to bound the sum
\begin{eqnarray*}
III=\left|\sum_{q'}\sum_{q: d(q)\geq r^{-s}s^n}\langle A_q \ast \mu_{\rho,0}^s, A_{q'} \ast \mu_{\rho,0}^s\rangle\right|.
\end{eqnarray*}
If $d(S)$ is the distance from $S$ to $q'$, then again
\begin{eqnarray*}
\left|\tilde\mu_{\rho,0}^s\ast\mu_{\rho,0}^s\ast \sum_{Q\subset S, \kappa(Q)=-s}\lambda_Qa_Q(x)\right|&\leq& C\sum_{Q\subset S^*, \kappa(Q)=-s} r^s|s|^n 2^{(5-d)\epsilon s}d(S)^{-2}\lambda_Q\\
&\leq& r^s|s|^n2^{(5-d)\epsilon s}\alpha|S|d(S)^{-2}
\end{eqnarray*}
using property (\ref{turned}) from the Whitney decomposition. Now, since the $S$ are disjoint, we can replace the sum over $S$ with $\int_{r^{-s}s^n}^5 t^{-2} t^{d-1} dt\leq C(1+s)$. Thus
\begin{eqnarray*}
III &\lesssim & \alpha \sum_Q \lambda_Q r^{-s}s^n2^{(5-d)\epsilon s}(1+s)
\end{eqnarray*}
which as above is good enough.

\begin{remark}
The reason that this argument works in $H^1$ and not in $L^1$ is that we used the cancellation of atoms to prove (\ref{pressure}).
\end{remark}

\section*{Acknowledgments}

This paper developed as a spinoff of a different project with A. Seeger; the author thanks him for the introduction to this problem and for many subsequent conversations on it.

\bibliography{BibLaVic12}{}
\bibliographystyle{plain}

\end{document}